\documentstyle[11pt]{article}
\def\Bbb R{{\rm \bf R}}
\def\proclaim#1{\vskip2mm{\bf #1}\em}
\def\endproclaim{\em \vskip2mm}
\def\tag#1{\eqno(#1)}
\def\gathered{\begin{array}{c}}
\def\endgathered{\end{array}}
\def\text{\mbox}

\begin{document}

\title {Remarks on the paper: Ikehata, M., Extraction formulae for an inverse boundary value problem for the equation $\nabla\cdot(\sigma-i\omega\epsilon)\nabla u=0$, Inverse Problems, {\bf 18}(2002), 1281-1290}
\author{Masaru IKEHATA\footnote{
Laboratory of Mathematics,
Graduate School of Adavanced Science and Engineering,
Hiroshima University, Higashhiroshima 739-8527, JAPAN}
\footnote{Emeritus Professor at Gunma University\,(2013 April 1st\,-\,)}}
\maketitle

\begin{abstract}
Some remarks on the jump condition appearing in Theorems 1.1 and 1.2 in the article
\cite{I1} and their implications are given.

\noindent
AMS: 35R30

\end{abstract}

\section{Recalling  Theorems 1.1 and 1.2 (\cite{I1})}

Let $\Omega$ be a bounded connected open subset
of $\Bbb R^n$, $n = 2,3$ with Lipschitz boundary. In what follows, unless otherwise stated, we assume
that $\sigma$, $\epsilon$ satisfy (1.1):

$$
\left\{
\begin{array}{c}
\displaystyle
\text{$\sigma$ and $\epsilon$ are $n\times n$ real symmetric matrix-valued functions on $\Omega$};\\
\displaystyle
\text{all components of  $\sigma$ and $\epsilon$ are essentially bounded functions on $\Omega$};\\
\displaystyle
\text{$\sigma$ is non-negative and $\epsilon$ is uniformly positive definite in $\Omega$.}
\end{array}
\right.
\tag {1.1}
$$

Given $f\in H^{1/2}(\partial D)$ there exists the unique weak
solution $u\in H^1(\Omega)$ of the Dirichlet problem
$$\begin{array}{c}
\displaystyle
\nabla\cdot(\sigma-i\omega\epsilon)\nabla u=0\,\,\text{in}\,\Omega,\\
\\
\displaystyle
u=f\,\,\text{on}\,\partial\Omega.
\end{array}
$$
Define the bounded linear functional $\Lambda_{\sigma,\epsilon}f$ on $H^{1/2}(\partial\Omega)$ by the formula
$$\displaystyle
<\Lambda_{\sigma,\epsilon}f,g>
=\int_{\Omega}(\sigma-i\omega\epsilon)\nabla u\cdot\nabla vdx
$$
where $g$ is an arbitrary element in $H^{1/2}(\partial\Omega)$ and $v\in H^1(\Omega)$ with $v=g$ on $\partial\Omega$.
The $\Lambda_{\sigma,\epsilon}$ is called the Dirichlet-to-Neumann map associated with the equation $\nabla\cdot(\sigma-i\omega\epsilon)\nabla u=0$.

Let $D$ be an open subset of $\Omega$ such that $\overline {D}\subset\Omega$.
Assume that $\sigma$, $\epsilon$ take the form

$$\displaystyle
\sigma(x)=\left\{
\begin{array}{ll}
\displaystyle
\sigma_0\,I_n, & \,\mbox{if $x\in\,\Omega\setminus D$,}\\
\\
\displaystyle
\sigma_0\,I_n+\alpha (x), & \,\mbox{if $x\in D$;}
\end{array}
\right.
\tag {1.2}
$$
$$\displaystyle
\epsilon(x)=\left\{
\begin{array}{ll}
\displaystyle
\epsilon_0\,I_n, & \,\mbox{if $x\in\,\Omega\setminus D$,}\\
\\
\displaystyle
\epsilon_0\,I_n+\beta (x), & \,\mbox{if $x\in D$,}
\end{array}
\right.
\tag {1.3}
$$
where $I_n$ denotes the $n\times n$-identity matrix, both $\sigma_0$ and $\epsilon_0$ are known constants satisfying
$$
\displaystyle
\sigma_0\ge 0;
\tag {1.4}
$$
$$
\displaystyle
\epsilon_0>0.
\tag {1.5}
$$
We assume that both $\alpha(x)$ and $\beta(x)$ together with $D$ are unknown and that
$(\sigma,\epsilon)$ has some kind of discontinuity across  $\partial D$ described below.

\subsection{A reduction procedure}

In this subsection we describe a simple reduction argument.
For $\sigma$ and $\epsilon$ given by (1.2) and (1.3), respectively define

$$\displaystyle
\tilde{\sigma}
=\frac{\sigma_0\sigma+\omega^2\epsilon_0\epsilon}
{\sigma_0^2+\omega^2\epsilon_0^2};
\tag {1.6}
$$
$$\displaystyle
\tilde{\epsilon}
=\frac{\sigma_0\epsilon-\epsilon_0\sigma}
{\sigma_0^2+\omega^2\epsilon_0^2}.
\tag {1.7}
$$
Then we have
$$\displaystyle
\sigma-i\omega\epsilon
=(\sigma_0-i\omega\epsilon_0)
(\tilde{\sigma}-i\omega\tilde{\epsilon}).
\tag {1.8}
$$
Note that $\tilde{\sigma}(x)=I_n$ and $\tilde{\epsilon}(x)=O_n$ for
$x\in\Omega\setminus D$.  From (1.1), (1.4)-(1.6) one knows that $\tilde{\sigma}$ is
uniformly positive definite in $\Omega$.  Then $\Lambda_{\tilde{\sigma},\tilde{\epsilon}}$
is still well defined and from (1.8) one has
$$\displaystyle
\Lambda_{\sigma,\epsilon}
=(\sigma_0-i\omega\epsilon_0)\Lambda_{\tilde{\sigma},\tilde{\epsilon}}.
\tag {1.9}
$$
Therefore, knowing $\Lambda_{\sigma,\epsilon}$ is equivalent to knowing $\Lambda_{\tilde{\sigma},\tilde{\epsilon}}$
through the relationship (1.9).
Moreover, from (1.6) and (1.7) we have
$$\displaystyle
\left(
\begin{array}{c}
\displaystyle
\tilde{\sigma}-I_n\\
\\
\displaystyle
\tilde{\epsilon}
\end{array}
\right)
=\frac{1}{\sigma_0^2+\omega^2\epsilon_0^2}
\left(
\begin{array}{lr}
\displaystyle
\sigma_0 & \omega^2\epsilon_0\\
\\
\displaystyle
-\epsilon_0 & \sigma_0
\end{array}
\right)
\left(
\begin{array}{c}
\displaystyle
\sigma-\sigma_0\\
\\
\displaystyle
\epsilon-\epsilon_0
\end{array}
\right).
$$
In other words, we have
$$\displaystyle
\tilde{\sigma}(x)=\left\{
\begin{array}{ll}
\displaystyle
I_n, & \,\mbox{if $x\in\,\Omega\setminus D$,}\\
\\
\displaystyle
I_n+a (x), & \,\mbox{if $x\in D$;}
\end{array}
\right.
\tag {1.10}
$$
$$\displaystyle
\tilde{\epsilon}(x)=\left\{
\begin{array}{ll}
\displaystyle
O_n, & \,\mbox{if $x\in\,\Omega\setminus D$,}\\
\\
\displaystyle
b (x), & \,\mbox{if $x\in D$,}
\end{array}
\right.
\tag {1.11}
$$
where $O_n$ denotes the $n\times n$-zero matrix, $a$ and $b$ are related to the $\alpha$ and $\beta$ in (1.2) and (1.3), respectively through the equations
$$\displaystyle
\left(
\begin{array}{c}
\displaystyle
a\\
\\
\displaystyle
b
\end{array}
\right)
=\frac{1}
{\sigma_0^2+\omega^2\epsilon_0^2}
\left(
\begin{array}{lr}
\displaystyle
\sigma_0 & \omega^2\epsilon_0\\
\\
\displaystyle
-\epsilon_0 & \sigma_0
\end{array}
\right)
\left(\begin{array}{c}
\displaystyle
\alpha\\
\\
\displaystyle
\beta
\end{array}
\right).
\tag {1.12}
$$

Note that in \cite{I1} hereafter we write $\tilde{\sigma}\rightarrow\sigma$ and $\tilde{\epsilon}\rightarrow \epsilon$
\footnote{This is the meaning of the words ``Hereafter we consider the reduced case unless
otherwise stated...'', line 9-10 up on page 1283 in \cite{I1}.}.
However, to avoid a confusion in explaning the meaning of the positive/negative jump condition
described below  we do not use 
such saving of the symbols.

\subsection{The enclosure method}

First recall notation and some definition.

We denote by $S^{n-1}$ the set of all unit vectors in $\Bbb R^n$.  The function $h_D$ defined by
the equation
$$\displaystyle
h_D(\vartheta)=\sup_{x\in D}x\cdot\vartheta,\,\,\vartheta\in S^{n-1}
$$
is called the support function of $D$.
For each $\vartheta\in S^{n-1}$ and a positive number $\delta$ set
$$\displaystyle
D_{\vartheta}(\delta)
=\{x\in D\,\vert\,h_D(\vartheta)-\delta<x\cdot\vartheta\le h_D(\vartheta)\}.
$$

{\bf\noindent Definition (The positive/negative jump condition for $\tilde{\sigma}$ given by (1.10)).}
Given $\vartheta\in S^{n-1}$ we say that
$\tilde{\sigma}$ given by (1.10) has a positive/negative jump on $\partial D$
from the direction $\vartheta$ if there exist constants $C_{\vartheta} > 0$ and $\delta_{\vartheta} > 0$
such that, for almost all $x\in D_{\vartheta}(\delta_{\vartheta})$ 
the lowest eigenvalue of $a(x)$/$-a(x)$ is greater than $C_{\vartheta}$.

$\quad$

{\bf\noindent The enclosure method.}  Assume that $\sigma$ and $\epsilon$ takes the form (1.2) and (1.3) with $\sigma_0$ and $\epsilon_0$
satisfying (1.4) and (1.5).

{\bf\noindent Preliminary Step 1.}  Given the original $\Lambda_{\sigma,\epsilon}$ compute $\Lambda_{\tilde{\sigma},\tilde{\epsilon}}$
via the formula (1.9), that is
$$\displaystyle
\Lambda_{\tilde{\sigma},\tilde{\epsilon}}
=\frac{1}{\sigma_0-i\omega\epsilon_0}
\Lambda_{\sigma,\epsilon},
$$
where $\tilde{\sigma}$ and $\tilde{\epsilon}$ are given by (1.10) and (1.11).

{\bf\noindent Preliminary Step 2.}  Given $\vartheta\in S^{n-1}$ take $\vartheta^{\perp}\in S^{n-1}$
perpendicular to $\vartheta$.  Given $\tau>0$ and $t\in\Bbb R$ compute the reduced indicator function
$$\displaystyle
I_{\vartheta,\vartheta^{\perp}}(\tau,t)
=e^{-2\tau t}
\text{Re}\,
<(\Lambda_{\tilde{\sigma},\tilde{\epsilon}}
-\Lambda_{I_n, O_n})
(e^{\tau x\cdot(\vartheta+i\vartheta^{\perp})}\vert_{\partial\Omega}),
\overline{e^{\tau x\cdot(\vartheta+i\vartheta^{\perp})}\vert_{\partial\Omega}}>.
\tag {1.13}
$$
In the theorems stated below we always assume that $\partial D$ is Lipschitz, $C^2$ in the case 
when $n=2, 3$, respectively\footnote{Making the regularity sharp is not the purpose of the article
\cite{I1}.  It is out of my interest since in that time this direction of the research, that is,
seeking a direct formula in inverse obstacle problems
was like a blue ocean!}.

\proclaim{\noindent Theorem 1.1.}
Assume that $\tilde{\sigma}$ has a positive jump on $\partial D$ from the direction $\vartheta$.
Then we have

 if $t>h_D(\vartheta)$, then 
$\displaystyle
\lim_{\tau\longrightarrow\infty}\vert I_{\vartheta,\vartheta^{\perp}}(\tau,t)\vert=0$;

 if $t<h_D(\vartheta)$, then 
$\displaystyle
\lim_{\tau\longrightarrow\infty}\vert I_{\vartheta,\vartheta^{\perp}}(\tau,t)\vert=\infty$;

 if $t=h_D(\vartheta)$, then 
$\displaystyle
\liminf_{\tau\longrightarrow\infty}\vert I_{\vartheta,\vartheta^{\perp}}(\tau,t)\vert>0$.

Moreover, the formula
$$\displaystyle
\lim_{\tau\longrightarrow\infty}
\frac{\log\vert I_{\vartheta,\vartheta^{\perp}}(\tau,t)\vert}
{2\tau}
=h_D(\vartheta)-t\,\,\,\,\,\,\forall t\in\Bbb R,
$$
is valid.

\endproclaim

Note that there is no restriction on $\omega$.
However, if $\tilde{\sigma}$ has a negative jump on $\partial D$ from
direction $\vartheta$, we do not know whether one can relax the condition (1.16) indicated below.

\proclaim{\noindent Theorem 1.2.}
Let $M>0$ and $m>0$ satisfy
$$
\displaystyle
\tilde{\sigma}(x)\xi\cdot\xi\ge m\vert\xi\vert^2\,\,\text{a.e.}\,x\in D\,\,\,\,\forall\xi\in\Bbb R^n
\tag {1.14}
$$
and
$$
\displaystyle
\vert b(x)\xi\vert\le M\vert\xi\vert\,\,\text{a.e.}\,x\in D\,\,\,\,\forall\xi\in\Bbb R^n.
\tag {1.15}
$$
Assume that $\tilde{\sigma}$ has a negative jump on $\partial D$ from the direction $\vartheta$ and that, for the constant
$C_{\vartheta}$ in the condition the frequency $\omega$ satisfies
$$\displaystyle
0\le\omega<\frac{\sqrt{mC_{\vartheta}}}{M}.
\tag {1.16}
$$
Then we have the same conclusion as that of Theorem 1.1.

\endproclaim

\section{The positive/negative jump condition and implications}

Now let us explain the meaning of the positive/negative jump condition across on $\partial D$ for $\tilde{\sigma}$ in terms of the original $\sigma$
and $\epsilon$ given by (1.2) and (1.3).

It follows from (1.12) that
$$\begin{array}{ll}
\displaystyle
a(x)=\frac{\sigma_0\,\alpha(x)+\omega^2\,\epsilon_0\beta(x)}{\sigma_0^2+\omega^2\epsilon_0^2}, & \text{a.e. $x\in D$}.
\end{array}
\tag {2.1}
$$
Since $\alpha(x)=\sigma(x)-\sigma_0\,I_n$ and $\beta(x)=\epsilon(x)-\epsilon_0\,I_n$ for $x\in D$ and $\sigma_0^2+\omega^2\epsilon_0^2>0$, we see that:

$\bullet$  $\tilde{\sigma}$ given by (1.10) has a positive jump on $\partial D$
from the direction $\vartheta$ if and only if there exist positive constants $C_{\vartheta}$ and $\delta_{\vartheta}$ such that,
for all $\xi\in\Bbb R^n$ and almost all $x\in D_{\vartheta}(\delta_{\vartheta})$
$$\displaystyle
\left\{\sigma_0\,(\sigma(x)-\sigma_0\,I_n)+\omega^2\,\epsilon_0(\epsilon(x)-\epsilon_0\,I_n)
\right\}
\xi\cdot\xi\ge (\sigma_0^2+\omega^2\epsilon_0^2)\,C_{\vartheta}\vert\xi\vert^2.
\tag {2.2}
$$

$\bullet$  $\tilde{\sigma}$ given by (1.10) has a negative jump on $\partial D$
from the direction $\vartheta$ if and only if there exist positive constants $C_{\vartheta}$ and $\delta_{\vartheta}$ such that,
for all $\xi\in\Bbb R^n$ and almost all $x\in D_{\vartheta}(\delta_{\vartheta})$
$$\displaystyle
-\left\{\sigma_0\,(\sigma(x)-\sigma_0\,I_n)+\omega^2\,\epsilon_0(\epsilon(x)-\epsilon_0\,I_n)
\right\}\xi\cdot\xi\ge (\sigma_0^2+\omega^2\epsilon_0^2)\,C_{\vartheta}\vert\xi\vert^2.
\tag {2.3}
$$

In Theorem 1.1 it is not important to know the concrete value of the constant $C_{\vartheta}$ in (2.2). Thus one can replace
the positive constant $(\sigma_0^2+\omega^2\epsilon_0^2)C_{\vartheta}$
with another positive constant.  Theorem 1.1 becomes

\proclaim{\noindent Theorem 1.1'.}  Let $\sigma_0$ and $\epsilon_0$ satisfy (1.4) and (1.5), respectively.
Assume that there exist positive constants $C_{\vartheta}'$ and $\delta_{\vartheta}$ such that, 
for all $\xi\in\Bbb R^n$ and almost all $x\in D_{\vartheta}(\delta_{\vartheta})$
$$\displaystyle
\left\{\sigma_0\,(\sigma(x)-\sigma_0\,I_n)+\omega^2\,\epsilon_0(\epsilon(x)-\epsilon_0\,I_n)
\right\}
\xi\cdot\xi\ge C_{\vartheta}'\,\vert\xi\vert^2.
\tag {2.4}
$$
Then we have

 if $t>h_D(\vartheta)$, then 
$\displaystyle
\lim_{\tau\longrightarrow\infty}\vert I_{\vartheta,\vartheta^{\perp}}(\tau,t)\vert=0$;

 if $t<h_D(\vartheta)$, then 
$\displaystyle
\lim_{\tau\longrightarrow\infty}\vert I_{\vartheta,\vartheta^{\perp}}(\tau,t)\vert=\infty$;

 if $t=h_D(\vartheta)$, then 
$\displaystyle
\liminf_{\tau\longrightarrow\infty}\vert I_{\vartheta,\vartheta^{\perp}}(\tau,t)\vert>0$.

Moreover, the formula
$$\displaystyle
\lim_{\tau\longrightarrow\infty}
\frac{\log\vert I_{\vartheta,\vartheta^{\perp}}(\tau,t)\vert}
{2\tau}
=h_D(\vartheta)-t\,\,\,\,\,\,\forall t\in\Bbb R,
$$
is valid.

\endproclaim

As a corollary, if the matrix valued function 
$\sigma_0\,(\sigma(x)-\sigma_0\,I_n)+\omega^2\,\epsilon_0(\epsilon(x)-\epsilon_0\,I_n)$ is uniformly positive definite on $D$, then all the formulae in Theorem 1.1' are valid for all directions $\vartheta$.

In contrast to Theorem 1.1, the value of the constant $C_{\vartheta}$ in (2.3) plays an important role as indicated in the constraint (1.16)
on $\omega$.  And also the constants $m$ and $M$ in (1.14) and (1.15).

Now consider (1.14), (1.15) and (1.16).
It follows from (1.6) that
the condition (1.14) has the expression
$$\displaystyle
\frac{\sigma_0\,\sigma(x)+\omega^2\,\epsilon_0\,\epsilon(x)}{\sigma_0^2+\omega^2\epsilon_0^2}\xi\cdot\xi
\ge m\vert\xi\vert^2\,\,\text{a.e.$x\in D$}\,\,\forall\xi\in\Bbb R^n.
$$
It follows from (1.12) that $b(x)$ for almost all $x\in D$ takes the form
$$\begin{array}{ll}
\displaystyle
b(x)
&
\displaystyle
=\frac{-\epsilon_0\alpha(x)+\sigma_0\beta(x)}{\sigma_0^2+\omega^2\epsilon_0^2}
\\
\\
\displaystyle
&
\displaystyle
=\frac{-\epsilon_0(\sigma(x)-\sigma_0\,I_n)+\sigma_0(\epsilon(x)-\epsilon_0\,I_n)}{\sigma_0^2+\omega^2\epsilon_0^2}\\
\\
\displaystyle
&
\displaystyle
=\frac{-\epsilon_0\,\sigma(x)+\sigma_0\,\epsilon(x)}{\sigma_0^2+\omega^2\epsilon_0^2}.
\end{array}
$$
Thus (1.15) becomes
$$\displaystyle
\left\Vert
\frac{\epsilon_0\,\sigma(x)-\sigma_0\,\epsilon(x)}{\sigma_0^2+\omega^2\epsilon_0^2}
\right\Vert\le M\,\,\text{a.e.$x\in D$},
$$
where $\Vert K\Vert=\sup_{\vert\xi\vert\le 1}\vert K\xi\vert$ for $n\times n$-matrix $K$.

Therefore Theorem 1.2 becomes

\proclaim{\noindent Theorem 1.2'.}  Let $\sigma_0$ and $\epsilon_0$ satisfy (1.4) and (1.5), respectively.
Assume that there exist positive constants $C_{\vartheta}$ and $\delta_{\vartheta}$ such that,
for all $\xi\in\Bbb R^n$ and almost all $x\in D_{\vartheta}(\delta_{\vartheta})$
$$\displaystyle
\frac{\sigma_0\,(\sigma(x)-\sigma_0\,I_n)+\omega^2\,\epsilon_0(\epsilon(x)-\epsilon_0\,I_n)}
{\sigma_0^2+\omega^2\epsilon_0^2}\xi\cdot\xi\le\,-C_{\vartheta}\vert\xi\vert^2.
\tag {2.5}
$$
Let $M>0$ and $m>0$ satisfy
$$\displaystyle
\frac{\sigma_0\,\sigma(x)+\omega^2\,\epsilon_0\,\epsilon(x)}{\sigma_0^2+\omega^2\epsilon_0^2}\xi\cdot\xi
\ge m\vert\xi\vert^2\,\,\text{a.e.$x\in D$}\,\,\forall\xi\in\Bbb R^n
\tag {2.6}
$$
and
$$\displaystyle
\left\Vert
\frac{\epsilon_0\,\sigma(x)-\sigma_0\,\epsilon(x)}{\sigma_0^2+\omega^2\epsilon_0^2}
\right\Vert\le M\,\,\text{a.e.$x\in D$}.
\tag {2.7}
$$
Let $\omega$ satisfy
$$\displaystyle
0\le\omega<\frac{\sqrt{m\,C_{\vartheta}}}{M}.
\tag {2.8}
$$
Then we have the same conclusion as that of Theorem 1.1.

\endproclaim

Some remarks on (2.5), (2.6) and (2.7) are in order.

Let both $\sigma_0$ and $\epsilon_0$ be positive instead of (1.4) and (1.5).  
Then one can rewrite
$$\displaystyle
\frac{\sigma_0\,(\sigma(x)-\sigma_0\,I_n)+\omega^2\,\epsilon_0(\epsilon(x)-\epsilon_0\,I_n)}
{\sigma_0^2+\omega^2\epsilon_0^2}
=P\left(\frac{\sigma(x)}{\sigma_0}-I_n\right)+Q\left(\frac{\epsilon(x)}{\epsilon_0}-I_n\right),
$$
where
$$\displaystyle
P=P(\sigma_0^2,\omega^2\sigma_0^2)=\frac{\sigma_0^2}{\sigma_0^2+\omega^2\epsilon_0^2},\,\,\,
Q=Q(\sigma_0^2,\omega^2\epsilon_0^2)=\frac{\omega^2\epsilon_0^2}{\sigma_0^2+\omega^2\epsilon_0^2}.
$$
Note that $P+Q=1$, $P>0$ and $Q>0$.  Thus the left-hand side on (2.5) is nothing but a {\it convex combination} of the {\it dimensionless quantities}
$\frac{\sigma(x)}{\sigma_0}-I_n$ and $\frac{\epsilon(x)}{\epsilon_0}-I_n$.

Thus, roughly speaking, in the case when the matrix $\frac{\sigma(x)}{\sigma_0}-I_n$ is positive/negative and matrix $\frac{\epsilon(x)}{\epsilon_0}-I_n$ negative/positive
the validilty of (2.5) for a positive constant $m$ is quite delicate.

(2.6) takes the form
$$\displaystyle
\left(P\,\frac{\sigma(x)}{\sigma_0}+Q\,\frac{\epsilon(x)}{\epsilon_0}\right)\xi\cdot\xi
\ge m\vert\xi\vert^2\,\,\text{a.e.$x\in D$}\,\,\forall\xi\in\Bbb R^n.
$$

(2.7) takes the form
$$\displaystyle
\frac{\sigma_0\epsilon_0}{\sigma_0^2+\omega^2\epsilon_0^2}
\left\Vert
\frac{\sigma(x)}{\sigma_0}-\frac{\epsilon(x)}{\epsilon_0}
\right\Vert
\le M\,\,\text{a.e.$x\in D$}.
\tag {2.9}
$$

Here let $R$ satisfy
$$\displaystyle
\text{ess. sup}_{x\in D}\left\Vert
\frac{\sigma(x)}{\sigma_0}-\frac{\epsilon(x)}{\epsilon_0}
\right\Vert
\le R.
$$
Then we have, for almost all $x\in D$
$$\displaystyle
\frac{\sigma_0\epsilon_0}{\sigma_0^2+\omega^2\epsilon_0^2}
\left\Vert
\frac{\sigma(x)}{\sigma_0}-\frac{\epsilon(x)}{\epsilon_0}
\right\Vert
\le
\frac{R}{2\omega}.
$$
Thus one can choose $M$ in (2.9) as
$$\displaystyle
M=\frac{R}{2\omega}.
$$
Then, (2.8) becomes
$$\displaystyle
0<\omega<\frac{2\omega\,\sqrt{m\,C_{\vartheta}}}{R}.
$$
This is equivalent to the inequality 
$$
\displaystyle
R<2\sqrt{m\,C_{\vartheta}}.
$$

Thus one gets a corollary of Theorem 1.2'.

\proclaim{\noindent Corollary 2.1.}  Let $\sigma_0>0$ and $\epsilon_0>0$.
Assume that there exist positive constants $C_{\vartheta}$ and $\delta_{\vartheta}$ such that,
for all $\xi\in\Bbb R^n$ and almost all $x\in D_{\vartheta}(\delta_{\vartheta})$
$$\displaystyle
\left\{P\left(\frac{\sigma(x)}{\sigma_0}-I_n\right)+Q\left(\frac{\epsilon(x)}{\epsilon_0}-I_n\right)\right\}
\xi\cdot\xi\le -C_{\vartheta}\vert\xi\vert^2.
$$
Let $m>0$ satisfy
$$\displaystyle
\left(P\,\frac{\sigma(x)}{\sigma_0}+Q\,\frac{\epsilon(x)}{\epsilon_0}\right)\xi\cdot\xi
\ge m\vert\xi\vert^2\,\,\text{a.e.$x\in D$}\,\,\forall\xi\in\Bbb R^n
$$
and assume that
$$\displaystyle
\text{ess. sup}_{x\in D}\left\Vert
\frac{\sigma(x)}{\sigma_0}-\frac{\epsilon(x)}{\epsilon_0}
\right\Vert
<2\sqrt{m\,C_{\vartheta}}.
\tag {2.10}
$$
Then we have the same conclusion as that of Theorem 1.1.

\endproclaim

Theorefore we have succeeded in dropping the constraint (2.8) on $\omega$ and instead introduced a kind of similarity condition (2.10)
on the relative conductivity and permittivity which are both {\it dimensionless}.

Of course, we have also a corollary of Theorem 1.1'.

\proclaim{\noindent Corollary 2.2.}  Let $\sigma_0>0$ and $\epsilon_0>0$.
Assume that there exist positive constants $C_{\vartheta}$ and $\delta_{\vartheta}$ such that,
for all $\xi\in\Bbb R^n$ and almost all $x\in D_{\vartheta}(\delta_{\vartheta})$
$$\displaystyle
\left\{P\left(\frac{\sigma(x)}{\sigma_0}-I_n\right)+Q\left(\frac{\epsilon(x)}{\epsilon_0}-I_n\right)\right\}
\xi\cdot\xi\ge C_{\vartheta}\vert\xi\vert^2.
$$
Then we have the same conclusion as that of Theorem 1.1.

\endproclaim

$$\quad$$

\centerline{{\bf Acknowledgment}}

The author was partially supported by Grant-in-Aid for
Scientific Research (C)(No. 17K05331) of Japan  Society for
the Promotion of Science.

$$\quad$$

\vskip1cm
\noindent
e-mail address

ikehata@hiroshima-u.ac.jp

\end{document}